\documentclass[11pt,oneside]{amsart}
\usepackage{amssymb} 
\newtheorem{thm}{Theorem}[section] 
\newtheorem{cor}[thm]{Corollary}
\newtheorem{lem}[thm]{Lemma} 
\newtheorem{prop}[thm]{Proposition}
\theoremstyle{definition}
\theoremstyle{remark}
\newtheorem{rem}[thm]{Remark} 
\theoremstyle{proof}

\numberwithin{equation}{section}

\newcommand{\norm}[1]{\left\Vert#1\right\Vert}
\newcommand{\abs}[1]{\left\vert#1\right\vert}
\newcommand{\set}[1]{\left\{#1\right\}}
\newcommand{\brac}[1]{\left(#1\right)}
\newcommand{\scalar}[1]{\left \langle #1 \right \rangle}
\newcommand{\sscalar}[1]{\langle #1 \rangle}
\newcommand{\Real}{\mathbb{R}}

\newcommand{\eps}{\varepsilon}
\newcommand{\To}{\longrightarrow}

\newcommand{\BP}{\mathcal{BP}}
\newcommand{\I}{\mathcal{I}}
\newcommand{\E}{\mathcal{E}}
\newcommand{\M}{\mathcal{M}}
\newcommand{\C}{\mathcal{C}}

\newcommand{\Vol}[1]{\textrm{Vol} (#1)}

\renewcommand{\S}{\mathcal{S}}

\newcommand{\ang}{\measuredangle}
\newcommand{\Proj}{\textrm{Proj}}
\newcommand{\nproj}{\overline{\Proj}}
\renewcommand{\Im}{\textrm{Im }}
\newcommand{\Ker}{\textrm{Ker }}

\begin{document}

\title{Generalized Intersection Bodies are Not Equivalent}


\author{Emanuel Milman}
\email{emanuel.milman@gmail.com}
\address{Department of Mathematics, The Weizmann Institute of Science.}
\thanks{Supported in part by BSF and ISF}

\begin{abstract}
In \cite{Koldobsky-I-equal-BP}, A. Koldobsky asked whether two types
of generalizations of the notion of an intersection-body, are in
fact equivalent. The structures of these two types of generalized
intersection-bodies have been studied in
\cite{EMilman-Generalized-Intersection-Bodies}, providing
substantial evidence for a positive answer to this question. The
purpose of this note is to construct a counter-example, which
provides a surprising negative answer to this question in a strong
sense. This implies the existence of non-trivial non-negative
functions in the range of the spherical Radon transform, and the
existence of non-trivial spaces which embed in $L_p$ for certain
negative values of $p$.
\end{abstract}


\maketitle


\section{Introduction}

Let $\Vol{L}$ denote the Lebesgue measure of a set $L \subset
\Real^n$ in its affine hull, and let $G(n,k)$ denote the Grassmann
manifold of $k$ dimensional subspaces of $\Real^n$. Let $D_n$ denote
the Euclidean unit ball, and $S^{n-1}$ the Euclidean sphere. All of
the bodies considered in this note will be assumed to be
centrally-symmetric star-bodies (even if the central-symmetry
assumption is omitted). A centrally-symmetric star-body $K$ is a
compact set with non-empty interior such that $K = -K$, $tK \subset
K$ for all $t \in [0,1]$, and such that its radial function
$\rho_K(\theta) = \max \{r \geq 0 \; | \; r\theta \in K \}$ for
$\theta \in S^{n-1}$ is an even continuous function on $S^{n-1}$.


This note concerns two generalizations of the notion of an
\emph{intersection body}, first introduced by E. Lutwak in
\cite{Lutwak-dual-mixed-volumes} (see also
\cite{Lutwak-intersection-bodies}). A star-body $K$ is said to be an
intersection body of a star-body $L$, if $\rho_K(\theta) = \Vol{L
\cap \theta^\perp}$ for every $\theta \in S^{n-1}$, where
$\theta^\perp$ is the hyperplane perpendicular to $\theta$.
$K$ is said to be an intersection body, if it is the limit in the
radial metric $d_r$ of intersection bodies $\{K_i\}$ of star-bodies
$\{L_i\}$, where $d_r(K_1,K_2) = \sup_{\theta \in S^{n-1}}
\abs{\rho_{K_1}(\theta)- \rho_{K_2}(\theta)}$. This is equivalent
(e.g. \cite{Lutwak-intersection-bodies}, \cite{Gardner-BP-5dim}) to
$\rho_K = R^*(d\mu)$, where $\mu$ is a non-negative Borel measure on
$S^{n-1}$, $R^*$ is the dual transform (as in (\ref{eq:duality111}))
to the Spherical Radon Transform $R:C(S^{n-1}) \rightarrow
C(S^{n-1})$, which is defined for $f\in C(S^{n-1})$ as:
\begin{equation} \label{eq:Radon111}
R(f)(\theta) = \int_{S^{n-1} \cap \theta^\perp} f(\xi)
d\sigma_{\theta}(\xi),
\end{equation}
where $\sigma_{\theta}$ the Haar probability measure on $S^{n-1}
\cap \theta^\perp$.

The notion of an intersection body has been shown to be
fundamentally connected to the Busemann-Petty Problem (first posed
in \cite{Busemann-Petty}), which asks whether two
centrally-symmetric convex bodies $K$ and $L$ in $\Real^n$
satisfying:
\begin{equation}\label{eq:Busemann-Petty}
\Vol{K\cap H} \leq \Vol{L \cap H} \; \; \forall H \in G(n,n-1)
\end{equation}
necessarily satisfy $\Vol{K} \leq \Vol{L}$. It was shown in
\cite{Lutwak-intersection-bodies}, \cite{Gardner-BP-5dim} that the
answer is equivalent to whether all centrally-symmetric convex
bodies in $\Real^n$ are intersection bodies, and in a series of
results (\cite{Larman-Rogers}, \cite{Ball-BP}, \cite{Bourgain-BP},
\cite{Giannopoulos-BP}, \cite{Papadimitrakis-BP},
\cite{Gardner-BP-5dim}, \cite{Gardner-BP-3dim},
\cite{Koldobsky-lp-intersection-bodies},
\cite{Zhang-Correction-4dim}, \cite{GKS}) that this is true for
$n\leq 4$, but false for $n \geq 5$.

\medskip

In \cite{Zhang-Gen-BP}, G. Zhang considered a generalization of the
Busemann-Petty problem, in which $G(n,n-1)$ in
(\ref{eq:Busemann-Petty}) is replaced by $G(n,n-k)$, where $k$ is
some integer between $1$ and $n-1$. Zhang showed that the
\emph{generalized $k$-codimensional Busemann-Petty problem} is also
naturally associated to a class of generalized intersection-bodies,
which will be referred to as \emph{$k$-Busemann-Petty bodies} (note
that these bodies are referred to as \emph{$n-k$-intersection
bodies} in \cite{Zhang-Gen-BP} and \emph{generalized
$k$-intersection bodies} in \cite{Koldobsky-I-equal-BP}), and that
the generalized $k$-codimensional problem is equivalent to whether
all centrally-symmetric convex bodies in $\Real^n$ are
$k$-Busemann-Petty bodies. It was shown in \cite{Bourgain-Zhang}
(see also \cite{Rubin-Zhang}), and later in
\cite{Koldobsky-I-equal-BP}, that the answer is negative for $k <
n-3$, but the cases $k=n-3$ and $k=n-2$ remain open (the case
$k=n-1$ is obviously true). Several partial answers to these cases
are known. It was shown in \cite{Zhang-Gen-BP} (see also
\cite{Rubin-Zhang}) that when $K$ is a centrally-symmetric convex
body of revolution then the answer is positive for the pair $K,L$
with $k=n-2,n-3$ and any star-body $L$. When $k=n-2$, it was shown
in \cite{Bourgain-Zhang} that the answer is positive if $L$ is a
Euclidean ball and $K$ is convex and sufficiently close to $L$. This
was extended in \cite{EMilman-Remark-On-BP-Problem}, where it was
shown that this is again true for $k=n-2$ and $k=n-3$, when $L$ is
an arbitrary star-body and $K$ is sufficiently close to a Euclidean
ball (but to an extent depending on its curvature). Several other
generalizations of the Busemann-Petty problem were treated in
\cite{Rubin-Zhang}, \cite{Zvavitch-BP-Arbitrary-Measures},
\cite{BP-Hyperbolic}, \cite{BP-Hyperbolic-low-dim}.

Before defining the class of $k$-Busemann-Petty bodies we shall need
to introduce the $m$-dimensional Spherical Radon Transform, acting
on spaces of continuous functions as follows:
\begin{eqnarray}
\nonumber R_m: C(S^{n-1}) \To C(G(n,m)) \\
\nonumber R_m(f) (E) = \int_{S^{n-1}\cap E} f(\theta)
d\sigma_E(\theta) ,
\end{eqnarray}
where $\sigma_E$ is the Haar probability measure on $S^{n-1} \cap
E$. It is well known (e.g. \cite{Helgason-Book}) that as an operator
on \emph{even} continuous functions, $R_m$ is injective. The dual
transform is defined on spaces of \emph{signed} Borel measures $\M$
by:
\begin{eqnarray}
\label{eq:duality111} & R_m^*: \M(G(n,m)) \To \M(S^{n-1}) & \\
\nonumber & \int_{S^{n-1}} f R_m^*(d\mu) = \int_{G(n,m)} R_m(f) d\mu
& \forall f \in C(S^{n-1}),
\end{eqnarray}
and for a measure $\mu$ with continuous density $g$, the transform
may be explicitly written in terms of $g$ (see \cite{Zhang-Gen-BP}):
\begin{eqnarray}
\nonumber R_m^* g (\theta) = \int_{\theta \in E \in G(n,m)} g(E)
d\nu_{m,\theta}(E) ,
\end{eqnarray}
where $\nu_{m,\theta}$ is the Haar probability measure on the
homogeneous space $\set{E \in G(n,m) \;|\; \theta \in E}$.

We shall say that a body $K$ is a $k$-Busemann-Petty body if
$\rho_K^k = R_{n-k}^*(d\mu)$ as measures in $\M(S^{n-1})$, where
$\mu$ is a non-negative Borel measure on $G(n,n-k)$. We shall denote
the class of such bodies by $\BP_k^n$. Choosing $k=1$, for which
$G(n,n-1)$ is isometric to $S^{n-1} / Z_2$ by mapping $H$ to
$S^{n-1}\cap H^\perp$, and noticing that $R$ is equivalent to
$R_{n-1}$ under this map, we see that $\BP_1^n$ is exactly the class
of intersection bodies.

\medskip

In \cite{Koldobsky-I-equal-BP}, a second generalization of the
notion of an intersection body was introduced by A. Koldobsky, who
studied a different analytic generalization of the Busemann-Petty
problem. Following \cite{Koldobsky-I-equal-BP}, a
centrally-symmetric star-body $K$ is said to be a $k$-intersection
body of a star-body $L$, if $\Vol{K \cap H^\perp} = \Vol{L \cap H}$
for every $H \in G(n,n-k)$. $K$ is said to be a $k$-intersection
body, if it is the limit in the radial metric of $k$-intersection
bodies $\{K_i\}$ of star-bodies $\{L_i\}$. We shall denote the class
of such bodies by $\I_k^n$. Again, choosing $k=1$, we see that
$\I_1^n$ is exactly the class of intersection bodies.

In \cite{Koldobsky-I-equal-BP}, Koldobsky considered the
relationship between these two types of generalizations, $\BP_k^n$
and $\I_k^n$, and proved that $\BP_k^n \subset \I_k^n$ (see also
\cite{EMilman-Generalized-Intersection-Bodies}). Koldobsky also
asked whether the opposite inclusion is equally true for all $k$
between 2 and $n-2$ (for 1 and $n-1$ this is true):

\smallskip
\noindent\textbf{Question (\cite{Koldobsky-I-equal-BP}):} Is it true
that $\BP_k^n = \I_k^n$ for $n \geq 4$ and $2\leq k \leq n-2$ ?\\

\vspace{-7pt} If this were true, as remarked by Koldobsky, a
positive answer to the generalized $k$-codimensional Busemann-Petty
problem for $k \geq n-3$ would follow, since for those values of $k$
any centrally-symmetric convex body in $\Real^n$ is known to be a
$k$-intersection body
(\cite{Koldobsky-correlation-inequality},\cite{Koldobsky-convex-is-n-3-intersection},
\cite{Koldobsky-I-equal-BP}).

\medskip

In \cite{EMilman-Generalized-Intersection-Bodies}, it was shown that
these two classes $\BP_k^n$ and $\I_k^n$ share many identical
structural properties, suggesting that it is indeed reasonable to
believe that $\BP_k^n = \I_k^n$. Using techniques from Integral
Geometry for the class $\BP_k^n$ and Fourier transform of
distributions techniques for the class $\I_k^n$, the following
structure Theorem was established (see
\cite{EMilman-Generalized-Intersection-Bodies} for an account of
particular cases which were known before). We define the
\emph{$k$-radial sum} of two star-bodies $L_1,L_2$ as the star-body
$L$ satisfying $\rho_L^{k} = \rho_{L_1}^{k} + \rho_{L_2}^{k}$.

\smallskip
\noindent\textbf{Structure Theorem
(\cite{EMilman-Generalized-Intersection-Bodies}) }\emph{Let $\C =
\I$ or $\C = \BP$ and $k,l=1,\ldots,n-1$. Then:
\begin{enumerate}
\item
$\C_k^n$ is closed under full-rank linear transformations,
$k$-radial sums and taking limit in the radial metric.
\item
$\C_1^n$ is the class of intersection-bodies in $\Real^n$, and
$\C_{n-1}^n$ is the class of all symmetric star-bodies in $\Real^n$.
\item
Let $K_1 \in \C_{k_1}^n$, $K_2 \in \C_{k_2}^n$ and $l = k_1 + k_2
\leq n-1$. Then the star-body $L$ defined by $\rho_L^{l} =
\rho_{K_1}^{k_1} \rho_{K_2}^{k_2}$ satisfies $L \in \C_{l}^n$. As
corollaries:
\begin{enumerate}
\item
$\C_{k_1}^n \cap \C_{k_2}^n \subset \C_{k_1+k_2}^n$ if $k_1+k_2 \leq
n-1$.
\item
$\C_{k}^n \subset \C_{l}^n$ if $k$ divides $l$.
\item
If $K \in \C_{k}^n$ then the star-body $L$ defined by $\rho_L =
\rho_K^{k/l}$ satisfies $L \in \C_{l}^n$ for $l \geq k$.
\end{enumerate}
\item
If $K \in \C_k^n$ then any $m$-dimensional central section $L$ of
$K$ (for $m>k$) satisfies $L \in \C_k^m$.
\end{enumerate}
}

Despite this and other evidence from
\cite{EMilman-Generalized-Intersection-Bodies} for a positive answer
to Koldobsky's question, we give the following negative answer. Let
$O(n)$ denote the orthogonal group on $\Real^n$. Recall that a
star-body $K$ is called a body of revolution if its radial function
$\rho_K \in C(S^{n-1})$ is invariant under the natural action of
$O(n-1)$ identified as some subgroup of $O(n)$.

\begin{thm} \label{thm:main1}
Let $n \geq 4$ and $2 \leq k \leq n-2$. Then there exists an
infinitely smooth centrally-symmetric body of revolution $K$ such
that $K \in \I_k^n$ but $K \notin \BP_k^n$.
\end{thm}

Note that Theorem \ref{thm:main1} does not imply a negative answer
to the unresolved cases $k=n-2,n-3$ (for $n\geq 5$) of the
generalized Busemann-Petty problem, which pertains to \emph{convex}
bodies. Indeed, the $K$ we construct cannot be a convex body in
those ranges of $k$, since as already mentioned, convex bodies of
revolution are known (\cite{Zhang-Gen-BP}, see also
\cite{Rubin-Zhang}) to belong to $\BP_{n-2}^n$ and $\BP_{n-3}^n$.
Theorem \ref{thm:main1} does however imply that if one wishes to
prove a positive answer to these unresolved cases by means of
comparing $k$-intersection bodies to $k$-Busemann-Petty bodies, it
is essential to restrict one's attention to convex bodies.

\smallskip

Let $I:C(G(n,k)) \rightarrow C(G(n,n-k))$ denote the operator
defined by $I(f)(E) = f(E^\perp)$ for all $E \in G(n,n-k)$. Let
$R_{n-k}(C(S^{n-1})) = \Im R_{n-k}$ denote the range of $R_{n-k}$.
As explained in Section \ref{sec:2}, Theorem \ref{thm:main1} can be
equivalently reformulated as follows:


\begin{thm} \label{thm:main2}
Let $n \geq 4$ and $2 \leq k \leq n-2$. Then there exists an
infinitely smooth function $g \in C(G(n,n-k))$ such that
$R_{n-k}^*(g) \geq 1$ and $(I \circ R_k)^*(g) \geq 1$ as functions
in $C(S^{n-1})$, but $g$ is not non-negative as a functional on
$R_{n-k}(C(S^{n-1}))$. In other words, there exists a non-negative
$h \in R_{n-k}(C(S^{n-1}))$ such that $\int_{G(n,n-k)} g(E) h(E)
d\eta_{n,n-k}(E) < 0$, where $\eta_{n,n-k}$ is the Haar probability
measure on $G(n,n-k)$. Moreover, both $g$ and $h$ can be chosen to
be invariant under the action of $O(n-1)$.
\end{thm}

\medskip

In \cite{EMilman-Generalized-Intersection-Bodies}, several
equivalent formulations to Koldobsky's question were obtained using
cone-duality and the Hahn-Banach Theorem for convex cones. Let
$C_{+}(S^{n-1})$ denote the cone of non-negative continuous
functions on the sphere, and let $R_{n-k}(C(S^{n-1}))_{+}$ denote
the cone of non-negative functions in the image of $R_{n-k}$. Let
$\overline{A}$ denote the closure of a set $A$ in the corresponding
normed space. Note that by the results from
\cite{EMilman-Generalized-Intersection-Bodies}, $\overline{\Im I
\circ R_k} = \overline{\Im R_{n-k}}$, and hence:
\[
\overline{R_{n-k}(C(S^{n-1}))_{+}} \supset
\overline{R_{n-k}(C_{+}(S^{n-1})) + I \circ R_k (C_{+}(S^{n-1}))}.\]
As formally verified in
\cite{EMilman-Generalized-Intersection-Bodies}, the dual formulation
to Theorem \ref{thm:main2} then reads:

\begin{thm} \label{thm:main3}
Let $n \geq 4$ and $2 \leq k \leq n-2$. Then:
\[
R_{n-k}(C(S^{n-1}))_{+} \setminus \overline{R_{n-k}(C_{+}(S^{n-1}))
+ I \circ R_k (C_{+}(S^{n-1}))} \neq \emptyset.
\]
In other words, there exists an (infinitely smooth) function $f \in
R_{n-k}(C(S^{n-1}))_{+}$ which can not be approximated (in
$C(G(n,n-k))$) by functions of the form $R_{n-k}(g) + I \circ
R_k(h)$ with $g,h \in C_{+}(S^{n-1})$.
\end{thm}

Other equivalent formulations using the language of Fourier
transforms of homogeneous distributions are given in section
\ref{sec:5}. We comment here that one such formulation pertains to
embeddings in $L_p$ for negative values of $p$. The definition of
embedding into such a space (for $-n<p<0$) was given by Koldobsky in
\cite{Koldobsky-I-equal-BP} by means of analytic continuation of the
usual definition for $p>0$. It is known (see Section \ref{sec:5})
that for $p \geq -1$ ($p \neq 0$) and for $-n<p\leq -n+1$, any
star-body $K$ such that $(\Real^n,\norm{\cdot}_K)$ embeds in $L_p$
can be generated by starting with the Euclidean ball $D_n$, applying
full-rank linear transformations, $(-p)$-radial sums and taking the
limit in the radial metric. Our results imply that $p=-1$ and
$p=-n+1$ are critical values for this property, and that this is no
longer true for $p=-k$, $2 \leq k \leq n-2$. In other words, there
exist ``non-trivial" $n$-dimensional spaces which embed in $L_{-k}$
for $2 \leq k \leq n-2$.

\medskip

The rest of this note is organized as follows. In Section
\ref{sec:2}, we provide some additional background which is required
to see why Theorem \ref{thm:main2} implies Theorem \ref{thm:main1}
and Theorem \ref{thm:main3}. In Section \ref{sec:3}, we develop
several formulas for the Spherical Radon Transform and its dual for
functions of revolution, i.e. functions invariant under the action
of $O(n-1)$. In Section \ref{sec:4}, we use these formulas to prove
Theorem \ref{thm:main2}, thereby constructing the desired
counter-example to Koldobsky's question. In Section \ref{sec:5}, we
give several additional equivalent formulations to Theorem
\ref{thm:main1} using the language of Fourier transforms of
homogeneous distributions.

\medskip

\noindent \textbf{Acknowledgments.} I would like to sincerely thank
my supervisor Prof. Gideon Schechtman for his guidance. I would also
like to thank Prof. Alexander Koldobsky for encouraging me to think
about bodies of revolution.

\section{Additional Background} \label{sec:2}

In this section, we summarize the relevant results needed for this
note. We also explain why Theorem \ref{thm:main1} and
\ref{thm:main3} follow from Theorem \ref{thm:main2}. We refer to
\cite{EMilman-Generalized-Intersection-Bodies} for more details.

For a star-body $K$ (not necessarily convex), we define its
Minkowski functional as $\norm{x}_K = \min \set{ t \geq 0 \; | \; x
\in t K}$. When $K$ is a centrally-symmetric convex body, this of
course coincides with the natural norm associated with it. Obviously
$\rho_K(\theta) = \norm{\theta}^{-1}_K$ for $\theta \in S^{n-1}$.

It was shown by Koldobsky in \cite{Koldobsky-I-equal-BP} that for a
star-body $K$ in $\Real^n$, $K \in \I_k^n$ iff $\norm{\cdot}_K^{-k}$
is a positive definite distribution on $\Real^n$, meaning that its
Fourier transform (as a distribution)
$(\norm{\cdot}_K^{-k})^{\wedge}$ is a non-negative Borel measure on
$\Real^n$. We refer the reader to Section \ref{sec:5} for more on
Fourier transforms of homogeneous distributions, as this will not be
of essence in the ensuing discussion. To translate this result to
the language of Radon transforms, it was shown in \cite[Corollary
4.2]{EMilman-Generalized-Intersection-Bodies} that for a infinitely
smooth star-body $K$ and a (signed) Borel measure $\mu \in
\M(G(n,n-k))$:
\begin{equation} \label{eq:Fourier-Radon}
\norm{\cdot}_K^{-k} = R^*_{n-k}(d\mu) \textrm{ iff }
(\norm{\cdot}_K^{-k})^\wedge = c(n,k) (I \circ R)^*_{k}(d\mu),
\end{equation}
where $c(n,k)$ is some positive constant and the equalities above
are interpreted as equalities between measures on $S^{n-1}$. Hence,
it follows (\cite[Lemma
5.3]{EMilman-Generalized-Intersection-Bodies}) that for an
infinitely smooth star-body $K$ in $\Real^n$, $K \in I_k^n$ iff
there exists a (possibly signed) Borel measure $\mu \in
\M(G(n,n-k))$, such that as measures $\rho_K^k = R_{n-k}^*(d\mu)
\geq 0$ and $(I \circ R_k)^*(d\mu) \geq 0$.

This should be compared with the definition of $k$-Busemann-Petty
bodies: $K \in \BP_k^n$ iff $\rho_K^k = R_{n-k}^*(d\mu)$ as measures
on $S^{n-1}$ for a \emph{non-negative} Borel measure $\mu \in
\M(G(n,n-k))$. Since for such a measure, $(I \circ R_k)^*(d\mu) \geq
0$, it follows that every infinitely smooth $k$-Busemann-Petty body
is also a $k$-intersection body, and this easily implies (see
\cite[Corollary 4.4]{EMilman-Generalized-Intersection-Bodies}) that
$\BP_k^n \subset \I_k^n$ in general, as first showed by Koldobsky in
\cite{Koldobsky-I-equal-BP}.

$R_{n-k}$ is known (e.g. \cite{Helgason-Book}) to be injective on
the space of \emph{even} functions in $C(S^{n-1})$, so by duality
$R_{n-k}^*$ is onto a dense subset of \emph{even} measures in
$\M(S^{n-1})$, which is known to include even measures with
infinitely smooth densities. However, it is important to note that
for $2\leq k \leq n-2$, the image of $R_{n-k}$ is not dense in
$C(G(n,n-k))$, and equivalently, $R^*_{n-k}$ has a non-trivial
kernel. The above implies that for any infinitely smooth star-body
$K$, we can find a measure $\mu$ such that $\rho_K^k =
R_{n-k}^*(d\mu)$, but if $2\leq k \leq n-2$ this measure will not
unique. Nevertheless, as a functional on $R_{n-k}(C(S^{n-1}))$, such
a measure $\mu$ is determined uniquely. The conclusion is that if we
need to determine whether $K \in \BP_k^n$ given a representation
$\rho_K^k = R_{n-k}^*(d\mu)$ for some measure $\mu \in
\M(G(n,n-k))$, a necessary and sufficient condition is that $\mu$ is
a non-negative functional on $R_{n-k}(C(S^{n-1}))$, i.e.
$\int_{G(n,n-k)} R_{n-k}(h)(E) d\mu(E) \geq 0$ for any $h \in
C(S^{n-1})$ such that $R_{n-k}(h) \geq 0$. Indeed, any non-negative
functional on $R_{n-k}(C(S^{n-1}))$ can be extended to a
non-negative functional on $C(G(n,n-k))$ by a version of the
Hahn-Banach Theorem (see the remarks before \cite[Lemma
5.2]{EMilman-Generalized-Intersection-Bodies} for more details).

The above discussion explains why Theorem \ref{thm:main1} is an
immediate consequence of Theorem \ref{thm:main2}. Given the
infinitely smooth function $g$ provided by Theorem \ref{thm:main2},
we define the centrally-symmetric star-body $K$ given by $\rho_K^k =
R_{n-k}^*(g)$. Note that this indeed defines a star-body since
$R_{n-k}^*(g)\geq 0$. In fact, $K$ is an infinitely smooth star-body
since it is known (e.g. \cite{Gelfand-Graev-Rosu}) that
$R_{n-k}^*(g)$ is an infinitely smooth function on $S^{n-1}$ if $g$
is infinitely smooth; and since $\rho_K^k = R_{n-k}^*(g) \geq 1$, it
follows that $\rho_K$ itself is infinitely smooth. In addition $K
\in \I_k^n$ since $(I \circ R_k)^*(g)\geq 0$. But since $g$ is not a
non-negative functional on $R_{n-k}(C(S^{n-1}))$, if follows that $K
\notin \BP_k^n$.

To explain why Theorem \ref{thm:main1} is equivalent to Theorem
\ref{thm:main3}, we recall another result from
\cite{EMilman-Generalized-Intersection-Bodies}. Denote $\M =
M(G(n,n-k))$ for short, and let:
\[
\M(\BP_k^n) = \set{\mu \in \M  ; \mu \text{ is a non-negative
functional on } R_{n-k}(C(S^{n-1}))},
\]
and:
\[
\M(\I_k^n) = \set{\mu \in \M ; R_{n-k}^*(d\mu) \geq 0 \text{ and }
(I \circ R_k)^*(d\mu) \geq 0}.
\]
It should already be clear from the above discussion that the
statement $\BP_k^n = \I_k^n$ is equivalent to the statement
$\M(\BP_k^n) = \M(\I_k^n)$. By the Hahn-Banach Theorem for convex
cones, it is not hard to see (\cite[Theorem
5.6]{EMilman-Generalized-Intersection-Bodies}) that the latter
statement is dual to:
\begin{equation} \label{eq:close-close}
\overline{R_{n-k}(C(S^{n-1}))_{+}} =
\overline{R_{n-k}(C_{+}(S^{n-1})) + I \circ R_k (C_{+}(S^{n-1}))}.
\end{equation}

As follows from (\ref{eq:Fourier-Radon}), $\Ker R^*_{n-k} = \Ker (I
\circ R_k)^*$, and therefore $\overline{\Im R_{n-k}} = \overline{\Im
I \circ R_k}$. This explains why the right-hand side of
(\ref{eq:close-close}) is always a subset of the left. Theorem
\ref{thm:main1} shows that it is a proper subset, implying Theorem
\ref{thm:main3}. Since this Theorem is attained using a convex
separation argument, we have no constructive way of finding the
function $f$ of the Theorem. Albeit, we can always find an
infinitely smooth $f$, since the subspace of infinitely smooth
functions in $R_{n-k}(C(S^{n-1}))$ is known to be dense in
$R_{n-k}(C(S^{n-1}))$, and hence in
$\overline{R_{n-k}(C(S^{n-1}))}$.

\section{Radon Transform for Functions of Revolution} \label{sec:3}

Fix $n \geq 3$ and $\xi_0 \in S^{n-1}$. We denote by
$O_{\xi_0}(n-1)$ the subgroup of $O(n)$ whose natural action on
$S^{n-1}$ leaves $\xi_0$ invariant, and by $C_{\xi_0}(S^{n-1})$ the
linear subspace of functions in $C_e(S^{n-1})$ invariant under
$O_{\xi_0}(n-1)$. Clearly $O_{\xi_0}(n-1)$ is isometric to $O(n-1)$.
We refer to members of $C_{\xi_0}(S^{n-1})$ as spherical functions
of revolution. For $\xi_1,\xi_2 \in S^{n-1}$, let
$\ang(\xi_1,\xi_2)$ denote the angle in $[0,\pi/2]$ between $\xi_1$
and $\xi_2$, i.e. $\cos \ang(\xi_1,\xi_2) =
\abs{\scalar{\xi_1,\xi_2}}$. We also denote $\ang(\xi_1,0) = \pi/2$.
Clearly $F \in C_{\xi_0}(S^{n-1})$ iff $F(\xi) = f(\ang(\xi,\xi_0))$
for $f \in C([0,\pi/2])$. In that case, we denote by $\tilde{f}\in
C([0,1])$ the function given by $\tilde{f}(\cos \theta) =
f(\theta)$, so $F(\xi) = \tilde{f}(\abs{\scalar{\xi,\xi_0}})$. We
denote the operator $T: C([0,\pi/2]) \rightarrow C([0,1])$ defined
by $T(f) = \tilde{f}$, for future reference. It is well known by
polar integration (e.g. \cite{VilenkinClassicBook}), that:
\begin{equation} \label{eq:polar-integration}
\int_{S^{n-1}} F(\xi) d\sigma_n(\xi) = c_n \int_{0}^{\pi/2}
f(\theta) \sin^{n-2}(\theta) d\theta = d_n \int_{0}^{1} \tilde{f}(t)
(1-t^2)^{\frac{n-3}{2}} dt, \
\end{equation}
where $\sigma_n$ is the Haar probability measure on $S^{n-1}$ and
$c_n$ is a constant whose value can be deduced by using $F \equiv f
\equiv \tilde{f} \equiv 1$.

For $E \in G(n,k)$ and $\xi \in S^{n-1}$, denote by $\Proj_E \xi$
the orthogonal projection of $\xi$ onto $E$, and by $\nproj_E \xi :=
\Proj_E \xi / |\Proj_E \xi|$ if $\Proj_E \xi \neq 0$, and $\nproj_E
\xi := 0$ otherwise. When $E = span(\xi_1)$ for $\xi_1 \in S^{n-1}$,
we may sometimes replace $E$ by $\xi_1$ in $\Proj_E$ and $\nproj_E$.
Denote by $\ang(\xi,E) = \ang(\xi,\nproj_E \xi)$ if $\nproj_E \xi
\neq 0$ and $\ang(\xi,E) = \pi/2$ otherwise.

Since the natural action of $O(n)$ on $C(G(n,k))$ and $C_e(S^{n-1})$
commutes with $R_k$, and since $O_{\xi_0}(n-1)$ acts transitively on
all $E \in G(n,k)$ such that $\ang(\xi_0,E)$ is fixed, it clearly
follows that if $F \in C_{\xi_0}(S^{n-1})$ then $R_k(F)(E)$ only
depends on $\ang(\xi_0,E)$. Hence, if $F(\xi) = f(\ang(\xi,\xi_0))$
for $f \in C[0,\pi/2]$, we denote (abusing notation) by $R_k(f) \in
C([0,\pi/2])$ the function given by $R_k(f)(\ang(\xi_0,E)) =
R_k(F)(E)$. Similarly, we define $\tilde{R}_k : C([0,1]) \rightarrow
C([0,1])$ as $\tilde{R}_k = T \circ R_k \circ T^{-1}$.

The following lemma was essentially stated in \cite{Zhang-Gen-BP}.
We provide a simple proof for completeness:

\begin{lem} \label{lem:R_k}
Let $f \in C[0,\pi/2]$ and $2 \leq k \leq n-1$. Then:
\[
R_k(f)(\phi) = c_k \int_0^{\pi/2} f(\cos^{-1}(\cos \phi \cos
\theta)) \sin^{k-2} \theta d\theta,
\]
where the value of $c_k$ is found by using $f \equiv 1$, in which
case $R_k(f) \equiv 1$.
\end{lem}
\begin{rem}
This lemma, together with the subsequent ones, extend to the case
$k=1$, if we properly interpret the (formally) diverging integral as
integration with respect to an appropriate delta-measure. Note also
that the value $c_k$ is consistent with the one used in
(\ref{eq:polar-integration}).
\end{rem}
\begin{proof}
Let $F\in C_{\xi_0}(S^{n-1})$ be given by $F(\xi) =
f(\ang(\xi,\xi_0))$. Let $E \in G(n,k)$ be such that $\ang(\xi_0,E)
= \phi$.
Hence, if $\xi_1 = \nproj_E \xi_0$ then $\ang(\xi_0,\xi_1) =
\phi$. For $\xi \in S^{n-1} \cap E$, since $\xi - \Proj_{\xi_1} \xi$
and $\xi_0 - \Proj_{\xi_1} \xi_0$ are orthogonal, it follows that
$\Proj_{\xi_0} \xi = \Proj_{\xi_0} (\Proj_{\xi_1} \xi)$. Hence $\cos
\ang(\xi,\xi_0) = \cos \ang(\xi,\xi_1) \cos \ang(\xi_1,\xi_0) = \cos
\ang(\xi,\xi_1) \cos \phi$. Since the function $F$ is even, a
standard polar integration formula then gives:
\begin{eqnarray}
\nonumber R_k(f)(\phi) &=& R_k(F)(E) = \int_{S^{n-1}\cap E} F(\xi)
d\mu_E(\xi) = \int_{S^{n-1}\cap E} f(\ang(\xi,\xi_0))
d\mu_E(\xi) \\
\nonumber &=& \int_{S^{n-1}\cap E} f(cos^{-1}(\cos \ang(\xi,\xi_1)
\cos \phi)) d\mu_E(\xi) \\
\nonumber &=& c_{n,k} \int_0^{\pi/2} f(\cos^{-1}(\cos \phi \cos
\theta)) \sin^{k-2} \theta d\theta.
\end{eqnarray}
\end{proof}

Performing the change of variables $t = \cos \theta$, $s = \cos
\phi$ above, we immediately have:

\begin{cor} \label{cor:tilde-R_k}
Let $\tilde{f} \in C[0,1]$ and $2 \leq k \leq n-1$. Then:
\[
\tilde{R}_k(\tilde{f})(s) = c_k \int_0^1 \tilde{f}(s t) (1-
t^2)^{\frac{k-3}{2}} dt,
\]
where the value of $c_k$ is the same as in Lemma \ref{lem:R_k}.
\end{cor}

Next, we introduce $C_{\xi_0}(G(n,k))$, the linear subspace of all
functions in $C(G(n,k))$ invariant under the action of
$O_{\xi_0}(n-1)$. We refer to members of $C_{\xi_0}(G(n,k))$ as
functions of revolution on the Grassmannian. As before, it is clear
that $G \in C_{\xi_0}(G(n,k))$ iff $G(E) = g(\ang(\xi_0,E))$ for $g
\in C([0,\pi/2])$. We have the following:

\begin{lem} \label{lem:bi-polar-integration}
Let $G \in C_{\xi_0}(G(n,k))$ such that $G(E) = g(\ang(\xi_0,E))$,
and let $\tilde{g} = T(g)$. Then:
\begin{eqnarray}
\nonumber \int_{G(n,k)} G(E) d\eta_{n,k}(E) & = & b_{n,k}
\int_{0}^{\pi/2} g(\phi) \sin^{n-k-1}\phi \cos^{k-1} \phi d\phi \\
\nonumber & = & b_{n,k} \int_{0}^{1} \tilde{g}(s)
(1-s^2)^{\frac{n-k-2}{2}} s^{k-1} ds,
\end{eqnarray}
Where $\eta_{n,k}$ is the Haar probability measure on $G(n,k)$, and
the value of $b_{n,k}$ may be deduce by using $G \equiv g \equiv
\tilde{g} \equiv 1$.
\end{lem}
\begin{proof}
Clearly:
\[
\int_{G(n,k)} G(E) d\eta_{n,k}(E) = \int_0^{\pi/2} g(\phi)
d\eta_{n,k}\set{E \in G(n,k) ; \ang(\xi_0,E) \leq \phi}.
\]
Since $\sigma_n$ and $\eta_{n,k}$ are rotation-invariant, it follows
that $\eta_{n,k}\set{E \in G(n,k) ; \ang(\xi_0,E) \leq \phi} =
\sigma_n\set{\xi \in S^{n-1} ; \ang(\xi,E_0) \leq \phi}$ for any
$E_0 \in G(n,k)$. Using bi-polar coordinates (e.g. \cite[Chapter
IX]{VilenkinClassicBook}), it is easy to see that:
\[
d\sigma_n\set{\xi \in S^{n-1} ; \ang(\xi,E_0) \leq \phi} = b_{n,k}
sin^{n-k-1} \phi cos^{k-1} \phi d\phi,
\]
for some $b_{n,k}$. This concludes the proof of the first equality
of the lemma, and the second one follows by the change of variables
$s = \cos(\phi)$.
\end{proof}

Next, we find an expression for the dual spherical Radon-Transform
of a function in $C_{\xi_0}(G(n,k))$. As before, it is clear that if
$F \in C_{\xi_0}(S^{n-1})$ then $R_k(F) \in C_{\xi_0}(G(n,k))$, and
that if $G \in C_{\xi_0}(G(n,k))$ then $R_k^*(G) \in
C_{\xi_0}(S^{n-1})$. If $G \in C_{\xi_0}(G(n,k))$ is given by $G(E)
= g(\ang(\xi_0,E))$, we denote by $R^*_k(g) \in C([0,\pi/2])$ the
function given by $R^*_k(g)(\ang(\xi,\xi_0)) = R^*_k(G)(\xi)$. As
usual, we define $\tilde{R}^*_k : C[0,1] \rightarrow C[0,1]$ by
$\tilde{R}^*_k = T \circ R^*_k \circ T^{-1}$. The standard duality
relation:
\[
\int_{S^{n-1}} R_k^*(G)(\xi) F(\xi) d\sigma_n(\xi) = \int_{G(n,k)}
G(E) R_k(F)(E) d\eta_{n,k}(E)
\]
is immediately translated using (\ref{eq:polar-integration}) and
Lemma \ref{lem:bi-polar-integration} into the following duality
relation between $\tilde{R}_k$ and $\tilde{R}^*_k$ on $C([0,1])$:

\begin{lem} \label{lem:revolution-duality}
Let $\tilde{f},\tilde{g} \in C([0,1])$ and $1 \leq k \leq n-1$.
Then:
\[
\int_0^1 \tilde{R}^*_k(\tilde{g})(t) \tilde{f}(t)
(1-t^2)^{\frac{n-3}{2}} dt = d_{n,k} \int_0^1 \tilde{g}(s)
\tilde{R}_k(\tilde{f})(s) (1-s^2)^{\frac{n-k-2}{2}} s^{k-1} ds
\]
where the value of $d_{n,k}$ is found by using $\tilde{f},\tilde{g}
\equiv 1$, in which case
$\tilde{R}_k(\tilde{f}),\tilde{R}^*_k(\tilde{g}) \equiv 1$.
\end{lem}

We can now deduce an expression for $\tilde{R}^*_k$:

\begin{lem} \label{lem:tilde-R^*_k}
Let $\tilde{g} \in C([0,1])$ and $2 \leq k \leq n-1$. Then:
\[
\tilde{R}^*_k(\tilde{g})(t) = e_{n,k} \int_0^1
\tilde{g}(\sqrt{1-s^2(1-t^2)}) (1-s^2)^{\frac{k-3}{2}} s^{n-k-1} ds,
\]
where the value of $e_{n,k}$ is found by using $\tilde{g} \equiv 1$,
in which case $\tilde{R}^*_k(\tilde{g}) \equiv 1$.
\end{lem}

\begin{proof}
We start with Lemma \ref{lem:revolution-duality} and use the formula
for $\tilde{R}_k$ given in Corollary \ref{cor:tilde-R_k}:
\begin{eqnarray}
\nonumber \int_0^1 \tilde{R}^*_k(\tilde{g})(t) \tilde{f}(t)
(1-t^2)^{\frac{n-3}{2}} dt = d_{n,k} \int_0^1 \tilde{g}(s)
\tilde{R}_k(\tilde{f})(s) (1-s^2)^{\frac{n-k-2}{2}} s^{k-1} ds \\
\nonumber = d_{n,k} c_k \int_0^1 \tilde{g}(s) \int_0^1 \tilde{f}(st)
(1-t^2)^{\frac{k-3}{2}} dt (1-s^2)^{\frac{n-k-2}{2}}
s^{k-1} ds \\
\nonumber = d_{n,k} c_k \int_0^1 \tilde{f}(v) \int_v^1
\tilde{g}(s)\brac{1-\frac{v^2}{s^2}}^{\frac{k-3}{2}}(1-s^2)^{\frac{n-k-2}{2}}
s^{k-2} ds dv.
\end{eqnarray}
Since this is true for any $\tilde{f} \in C([0,1])$, setting
$e_{n,k} = d_{n,k} c_k$, we conclude that:
\[
\tilde{R}^*_k(\tilde{g})(t) = e_{n,k} (1-t^2)^{-\frac{n-3}{2}}
\int_t^1
\tilde{g}(s)\brac{1-\frac{t^2}{s^2}}^{\frac{k-3}{2}}(1-s^2)^{\frac{n-k-2}{2}}
s^{k-2} ds.
\]
By the change of variable $s = \sqrt{1 - (s')^2(1-t^2)}$, one easily
checks that the assertion of the lemma is obtained.
\end{proof}

We now recall the definition of the ``perp" operator $I$ from the
Introduction, and extend it to the context of functions of
revolution. For every $k=1,\ldots,n-1$, we define $I: C(G(n,k))
\rightarrow C(G(n,n-k))$ as $I(f)(E) = f(E^\perp)$ for all $E \in
G(n,n-k)$, without specifying the index $k$. $I$ is obviously
self-adjoint:
\[
\int_{G(n,n-k)} I(F)(H) G(H) d\eta_{n-k}(H) = \int_{G(n,k)} F(E)
I(G)(E) d\eta_k(E),
\]
for all $F \in C(G(n,k))$ and $G \in C(G(n,n-k))$, where $\eta_m$
denotes the Haar probability measure on $G(n,m)$.

Since $\ang(\xi_0,E) = \pi/2 - \ang(\xi_0,E^\perp)$, it is clear
that for $G \in C_{\xi_0}(G(n,k))$ such that $G(E) =
g(\ang(\xi_0,E))$ for every $E \in G(n,k)$, $I(G)(H) = g(\pi/2 -
\ang(\xi_0,H))$ for every $H \in G(n,n-k)$. We therefore define $I :
C([0,\pi/2]) \rightarrow C([0,\pi/2])$ as $I(g)(\phi) =
g(\pi/2-\phi)$. Similarly, for $\tilde{g} \in C([0,1])$, we define
$I(\tilde{g})(s) = \tilde{g}(\sqrt{1-s^2})$. Clearly, if $G(E) =
\tilde{g}(cos(\ang(\xi_0,E)))$ then $I(G)(H) =
I(\tilde{g})(cos(\ang(\xi_0,H)))$. Hence in both cases $I$ must be
self-adjoint, and this can be also verified directly. As an
immediate corollary of \ref{lem:tilde-R^*_k}, we have:

\begin{cor}
Let $\tilde{g} \in C([0,1])$ and $2 \leq k \leq n-1$. Then:
\[
(I \circ \tilde{R}_k)^*(\tilde{g})(t) = e_{n,k} \int_0^1 \tilde{g}(s
\sqrt{1-t^2}) (1-s^2)^{\frac{k-3}{2}} s^{n-k-1} ds,
\]
where the value of $e_{n,k}$ is the same as in Lemma
\ref{lem:tilde-R^*_k}.
\end{cor}

We are now ready to construct the counter-example to Koldobsky's
question, as described in the next section.

\section{The Construction} \label{sec:4}

The main step in the proof of Theorem \ref{thm:main2}, is the
following:

\begin{prop} \label{prop:main}
For any $n \geq 4$, $2 \leq k \leq n-2$ and $s_0 \in (0,1)$, there
exists an infinitely smooth function $\tilde{g} \in C([0,1])$ such
that:
\begin{enumerate}
\item
For all $t \in [0,1]$:
\[
\tilde{R}_{n-k}^*(\tilde{g})(t) = e_{n,n-k} \int_0^1
\tilde{g}(\sqrt{1-s^2(1-t^2)})(1-s^2)^{\frac{n-k-3}{2}} s^{k-1} ds
\geq 1.
\]
\item
For all $t \in [0,1]$:
\[
(I \circ \tilde{R}_{k})^*(\tilde{g})(t) = e_{n,k} \int_0^1
\tilde{g}(s\sqrt{1-t^2})(1-s^2)^{\frac{k-3}{2}} s^{n-k-1} ds \geq 1.
\]
\item
$\tilde{g}(s_0) = -1.$
\end{enumerate}
\end{prop}
\begin{proof}
The proof is straightforward. We provide the details nevertheless.
Let $\eps>0$ be such that $[s_0-2\eps,s_0+2\eps] \subset (0,1)$. Let
$T_t,T'_t \in C([0,1])$ be defined by $T_t(s) = \sqrt{1-s^2(1-t^2)}$
and $T'_t(s) =s \sqrt{1-t^2}$, and let $\lambda$ denote the Lebesgue
measure on $\Real$. It is elementary to check that the maximum of
$\lambda\set{T_t^{-1}[s_0-2\eps,s_0+2\eps]}$ over $t\in[0,1]$ is
attained at $t = s_0-2\eps$, in which case it is equal to:
\[
\delta_1 := \max_{t \in [0,1]}
\lambda\set{T_t^{-1}[s_0-2\eps,s_0+2\eps]} = 1 -
\sqrt{\frac{1-(s_0+2\eps)^2}{1-(s_0-2\eps)^2}} < 1.
\]
An analogous computation shows that the maximum of
$\lambda\set{(T'_t)^{-1}[s_0-2\eps,s_0+2\eps]}$ over $t\in[0,1]$ is
attained at $t = \sqrt{1-(s_0+2\eps)^2}$, in which case it is equal
to:
\[
\delta_2 := \max_{t \in [0,1]}
\lambda\set{(T'_t)^{-1}[s_0-2\eps,s_0+2\eps]} =
\frac{4\eps}{s_0+2\eps} < 1.
\]
Set $\delta := \max(\delta_1,\delta_2) < 1$. Now denote by
$\mu_{n,m}$ the measure $e_{n,m} (1-s^2)^{\frac{m-3}{2}} s^{n-m-1}
ds$ on $[0,1]$, for $2 \leq m \leq n-2$. These are probability
measures, as witnessed by using $\tilde{g} \equiv 1$ in Lemma
\ref{lem:tilde-R^*_k}, in which case $\tilde{R}_k^*(\tilde{g})
\equiv 1$. Since their densities (with respect to $\lambda$) are
absolutely continuous and do not vanish on $(0,1)$, a compactness
argument shows that (fixing $n$):
\[
\gamma := \sup_{v \in [0,1],2\leq m \leq n-2}
\mu_{n,m}([v,v+\delta]) < 1.
\]
Set $\gamma^* = \frac{1+\gamma}{1-\gamma}$. We conclude by
constructing $\tilde{g}$ as follows. Set $\tilde{g}(s) = -1$ for $s
\in [s_0-\eps,s_0+\eps]$, $\tilde{g}(s) = \gamma^*$ for $s \in [0,1]
\setminus [s_0-2\eps,s_0+2\eps]$, and for $s \in
[s_0-2\eps,s_0+2\eps] \setminus [s_0-\eps,s_0+\eps]$ set
$\tilde{g}(s) \in [-1,\gamma^*]$ so that the resulting function
$\tilde{g} \in C[0,1]$ is in fact infinitely smooth (using standard
methods). Alternatively, we could simply define $\tilde{g}(s) =
(\gamma^*+1)(\frac{s-s_0}{2\eps})^2 - 1$ on $[0,1]$. Setting:
\[
\beta_1(t) := \mu_{n,n-k}\set{s\in[0,1]; T_t(s) \in
[s_0-2\eps,s_0+2\eps]},
\]
the definition of $\gamma$ and $\delta$ imply that $\beta_1(t) \leq
\gamma$ for all $t \in [0,1]$, hence:
\[
\int_0^1 \tilde{g}(\sqrt{1-s^2(1-t^2)}) d\mu_{n,n-k}(s) \geq
\gamma^* (1-\beta_1(t)) - \beta_1(t) \geq 1
\]
for all $t \in [0,1]$. Similarly, setting:
\[
\beta_2(t) := \mu_{n,k}\set{s\in[0,1]; T'_t(s) \in
[s_0-2\eps,s_0+2\eps]},
\]
we have $\beta_2(t) \leq \gamma$ for all $t \in [0,1]$, and:
\[
\int_0^1 \tilde{g}(s\sqrt{1-t^2}) d\mu_{n,k}(s) \geq \gamma^*
(1-\beta_2(t)) - \beta_2(t) \geq 1
\]
for all $t \in [0,1]$. This concludes the proof.
\end{proof}
\begin{rem}
Note that for $k=1$ and $k=n-1$ the above reasoning fails, as the
measure $\mu_{n,1}$ is a singular measure.
\end{rem}

\begin{rem}
Note also that the function $\tilde{g}$ we have constructed in fact
satisfies the claims (1) and (2) for \emph{all} values of $k$ in the
range $2 \leq k \leq n-2$.
\end{rem}

We can now almost conclude the proof of Theorem \ref{thm:main2}. We
still need one last observation, since a-priori, the fact that
$\tilde{g}(s_0) < 0$ does not guarantee that the function $G \in
C(G(n,n-k))$ defined as $G(E) = \tilde{g}(\cos(\ang(\xi_0,E)))$, is
not a non-negative functional on $R_{n-k}(C(S^{n-1}))$. This is
resolved by the following:

\begin{lem} \label{lem:polynomials}
The polynomials on $[0,1]$ are in the range of
$\tilde{R}_{n-k}(C([0,1]))$.
\end{lem}
\begin{proof}
This is immediate by Corollary \ref{cor:tilde-R_k}, because if
$\tilde{p}(t) = t^m$ ($m\geq 0$), then:
\[
\tilde{R}_k(\tilde{p})(s) = c_k \int_0^1 \tilde{p}(s t) (1-
t^2)^{\frac{k-3}{2}} dt = d_{k,m} s^m,
\]
with $d_{k,m} > 0$. Hence polynomials are mapped to polynomials by
$\tilde{R}_{n-k}$, and any polynomial in the range may be obtained.
\end{proof}

By the Weierstrass approximation theorem, if follows that:
\begin{cor} \label{cor:dense-range}
The range of $\tilde{R}_{n-k}$ is dense in $C([0,1])$.
\end{cor}

We can now turn to the proof of Theorem \ref{thm:main2}.

\begin{proof}[Proof of Theorem \ref{thm:main2}]
Let $\tilde{g} \in C[0,1]$ be the infinitely smooth function
constructed in Proposition \ref{prop:main}, with, say $s_0=1/2$. Fix
some $\xi_0 \in S^{n-1}$, and let $G \in C_{\xi_0}(G(n,n-k))$ be
defined by $G(E) = \tilde{g}(\cos(\ang(\xi_0,E)))$ for every $E \in
G(n,n-k)$. Since the functions $\tilde{g}$, $\cos$ and
$\ang(\xi,\cdot)$ are infinitely smooth on their corresponding
domains, so is their composition, hence $G$ is infinitely smooth on
$G(n,n-k)$. By the construction of $\tilde{g}$ and the compatibility
of $R^*_{n-k}$ and $(I \circ R_k)^*$ with $\tilde{R}^*_{n-k}$ and
$(I \circ \tilde{R}_k)^*$, respectively, it follows that
$R^*_{n-k}(G) = \tilde{R}^*_{n-k}(\tilde{g}) \geq 1$ and $(I \circ
R)^*_{k}(G) = (I \circ \tilde{R})^*_{k}(\tilde{g}) \geq 1$. It
remains to show that $G$ is not a non-negative functional on
$R_{n-k}(C(S^{n-1}))$. Let $H \in C_{\xi_0}(S^{n-1})$ be such that
$H(\xi) = \tilde{h}(\cos(\ang(\xi_0,\xi))$ for some $\tilde{h} \in
C([0,1])$. Then by Lemma \ref{lem:bi-polar-integration}:
\begin{equation} \label{eq:stam}
\int_{G(n,n-k)} G(E) R_{n-k}(H)(E) d\eta_{n,k} = \int_0^1
\tilde{g}(s) \tilde{R}_{n-k}(\tilde{h})(s) (1-s^2)^{\frac{n-k-2}{2}}
s^{k-1} ds.
\end{equation}
Since $\tilde{g}(s)(1-s^2)^{\frac{n-k-2}{2}} s^{k-1}$ is a
continuous function on $[0,1]$ whose value at $s_0$ is negative, by
Corollary \ref{cor:dense-range} we can find a function $\tilde{h}
\in C([0,1])$ such that the integral in (\ref{eq:stam}) is negative.
This concludes the proof.
\end{proof}

\section{Additional formulations} \label{sec:5}

In this section, we give several additional equivalent formulations
to the main result of this note, using the language of Fourier
transforms of homogeneous distributions (we refer the reader to
\cite{Koldobsky-Book} for more on this subject).

We denote by $\S(\Real^n)$ the space of rapidly decreasing
infinitely differentiable test functions in $\Real^n$, and by
$\S'(\Real^n)$ the space of distributions over $\S(\Real^n)$. The
Fourier transform $\hat{f}$ of a distribution $f \in \S'(\Real^n)$
is defined by $\sscalar{\hat{f},\phi} = \sscalar{f,\hat{\phi}}$ for
every test function $\phi$, where $\hat{\phi}(y) = \int \phi(x)
\exp(-i\sscalar{x,y}) dx$. A distribution $f$ is called homogeneous
of degree $p \in \Real$ if $\scalar{f,\phi(\cdot/t)} = \abs{t}^{n+p}
\scalar{f,\phi}$ for every $t>0$, and it is called even if the same
is true for $t=-1$. An even distribution $f$ always satisfies
$(\hat{f})^\wedge = (2\pi)^n f$. The Fourier transform of an even
homogeneous distribution of degree $p$ is an even homogeneous
distribution of degree $-n-p$. A distribution $f$ is called positive
if $\sscalar{f,\phi} \geq 0$ for every $\phi \geq 0$, implying that
$f$ is necessarily a non-negative Borel measure on $\Real^n$. We use
Schwartz's generalization of Bochner's Theorem
(\cite{Gelfand-Shilov}) as a definition, and call a homogeneous
distribution positive-definite if its Fourier transform is a
positive distribution.

\smallskip

The following characterization was given by Koldobsky in
\cite{Koldobsky-I-equal-BP}:

\begin{thm}[Koldobsky] \label{thm:I-char}
The following are equivalent for a centrally-symmetric star-body $K$
in $\Real^n$:
\begin{enumerate}
\item $K$ is a $k$-intersection body.
\item $\norm{\cdot}_K^{-k}$ is a positive definite distribution on
$\Real^n$, meaning that its Fourier-transform
$(\norm{\cdot}_K^{-k})^{\wedge}$ is a non-negative Borel measure on
$\Real^n$.
\item
The space $(\Real^n,\norm{\cdot}_K)$ embeds in $L_{-k}$.
\end{enumerate}
\end{thm}

For completeness, we give the definition of embedding in $L_{-k}$.
For $p>-1$ (and $p \neq 0$, the case $p=0$ requires passing to the
limit), it is well known (e.g. \cite{Koldobsky-I-equal-BP}) that
$(\Real^n,\norm{\cdot})$ embeds in $L_p$ iff:
\begin{equation} \label{eq:SL_p^n}
\norm{x}^p = \int_{S^{n-1}} \abs{\scalar{x,\theta}}^p d\mu(\theta),
\end{equation}
for some $\mu \in \M_+(S^{n-1})$, the cone of non-negative Borel
measures on $S^{n-1}$. Unfortunately, this characterization breaks
down at $p=-1$ since the above integral no longer converges.
However, Koldobsky showed that it is possible to regularize this
integral by using Fourier-transforms of distributions, and gave the
following definition: $(\Real^n,\norm{\cdot})$ embeds in $L_{-p}$
for $0<p<n$ iff there exists a measure $\mu \in \M_+(S^{n-1})$ such
that for any even test-function $\phi$:
\begin{equation} \label{eq:L_-k-characterization}
\int_{\Real^n} \norm{x}^{-p} \phi(x) dx = \int_{S^{n-1}}
\int_0^\infty t^{p-1} \hat{\phi}(t \theta) dt d\mu(\theta).
\end{equation}

In addition to the characterization (3) in Theorem \ref{thm:I-char}
of $\I_k^n$ as the class of unit-balls of subspaces of \emph{scalar}
$L_{-k}$ spaces, a functional analytic characterization of $\BP_k^n$
as the class of unit-balls of subspaces of certain
\emph{vector-valued} $L_{-k}$ spaces
was given in \cite{Koldobsky-I-equal-BP}. To explain this better, we
state the definition given by Koldobsky: $(\Real^n,\norm{\cdot})$
embeds in $L_{-p}(\Real^k)$ for $0<p<n$ iff there exists a measure
$\mu \in \M_+(\Real^{nk})$ such that for any even test-function
$\phi$:
\begin{equation} \label{eq:L_-k-vector-characterization}
\int_{\Real^n} \norm{x}^{-p} \phi(x) dx = \int_{\Real^{nk}}
\int_{\Real^k} \norm{v}_2^{p-k} \hat{\phi}\brac{\sum_{i=1}^k v_i
\xi_i} dv d\mu(\xi).
\end{equation}
For $k=1$ it is easy to see that this coincides with the definition
of embedding in $L_{-p}$. Using this definition, it was shown in
\cite{Koldobsky-I-equal-BP} that $K \in \BP_k^n$ iff
$(\Real^n,\norm{\cdot}_K)$ embeds in $L_{-k}(\Real^k)$.

For $p>0$, it is known that every separable vector valued $L_p$
space is isometric to a subspace of a scalar $L_p$ space and
vice-versa. Translating Theorem \ref{thm:main1} into the language of
$L_p$ spaces, we see that this is no longer true when $p = -k$, $2
\leq k \leq n-2$:

\begin{cor}
Let $n \geq 4$ and $2 \leq k \leq n-2$. Then there exists an
infinitely smooth centrally-symmetric body of revolution $K$ such
that $(\Real^n,\norm{\cdot}_K)$ embeds in $L_{-k}$ but does not
embed in $L_{-k}(\Real^k)$.
\end{cor}

\medskip

Next, we describe another property of $L_p$ spaces which breaks down
when passing the critical value of $p=-1$. Let us denote the set of
all star-bodies $K$ in $\Real^n$ for which
$(\Real^n,\norm{\cdot}_K)$ embeds in $L_p$ ($p \neq 0$) by $SL_p^n$.
For $p\neq 0$, let the $p$-norm sum of two bodies $L_1,L_2$ be
defined as the body $L$ satisfying $\norm{\cdot}_L^p =
\norm{\cdot}_{L_1}^p + \norm{\cdot}_{L_2}^p$. Obviously, the
$p$-norm sum coincides with the $(-p)$-radial sum, defined in the
introduction (before the Structure Theorem). We will denote by
$D_p^n$, the class of bodies created from the Euclidean ball $D_n$
by applying full-rank linear-transformations, $p$-norm sums, and
taking the limit in the radial metric. Using the characterization in
(\ref{eq:SL_p^n}), it is easy to show (e.g. \cite[Theorem
6.13]{Grinberg-Zhang}) that for $p>-1$ ($p \neq 0$), $SL_p^n =
D_p^n$. In order to understand what happens when $p \leq -1$, we
turn to the following characterization of $\BP_k^n$, first proved by
Goodey and Weil in \cite{Goodey-Weil} for intersection-bodies (the
case $k=1$), and extended to general $k$ by Grinberg and Zhang in
\cite{Grinberg-Zhang}:

\begin{thm}[Grinberg and Zhang] \label{thm:G&Z}
A star-body $K$ in $\Real^n$ is a $k$-Busemann-Petty body ($1 \leq k
\leq n-1$) iff it is the limit of $\set{K_i}$ in the radial metric,
where each $K_i$ is a finite $k$-radial sum of ellipsoids
$\set{\E^i_j}$ in $\Real^n$ having non-empty interior:
\[
\rho^k_{K_i} = \rho^k_{\E^i_1} + \ldots + \rho^k_{\E^i_{m_i}}.
\]
\end{thm}

In other words, Theorem \ref{thm:G&Z} states that $D_{-k}^n =
\BP_k^n$ for $k=1,\ldots,n-1$. Recall that $\I_1^n = \BP_1^n$ is the
class of all intersection-bodies in $\Real^n$ and $\I_{n-1}^n =
\BP_{n-1}^n$ is the class of all centrally-symmetric star-bodies in
$\Real^n$ (this is clear from the definitions, see also the
Structure Theorem from the introduction). Since $\I_k^n = SL_{-k}^n$
by characterization (3) of Theorem \ref{thm:I-char}, we see that
$SL_{-k}^n = D_{-k}^n$ for $k=1$ and $k=n-1$. However, Theorem
\ref{thm:main1} implies that this is no longer true for $2 \leq k
\leq n-2$:

\begin{cor} \label{cor:last-2}
Let $n \geq 4$ and $2 \leq k \leq n-2$. Then $SL_{-k}^n \setminus
D_{-k}^n \neq \emptyset$.
\end{cor}

Note that since $\BP_k^n \subset \I_k^n$, it is always true that
$D_{-k}^n \subset SL_{-k}^n$ (in fact, this is straightforward to
check directly, implying that $\BP_k^n \subset \I_k^n$ by using
Theorems \ref{thm:I-char} and \ref{thm:G&Z}). In some sense, the
members of $D_{-k}^n$ are the ``trivial" elements of $SL_{-k}^n$,
since obviously $D_n \in SL_{-k}^n$, and $SL_{-k}^n$ is closed under
taking full-rank linear transformations, $(-k)$-norm sums and and
limit in the radial-metric. Corollary \ref{cor:last-2} therefore
says that there are also ``non-trivial" elements in $SL_{-k}^n$, for
$2 \leq k \leq n-2$.

\medskip

We conclude by translating Corollary \ref{cor:last-2} into the
language of Fourier transforms of homogeneous distributions. Given
an even $f \in C(S^{n-1})$, we denote by $E_p(f)$ its homogeneous
extension of degree $p$ onto $\Real^n$ (formally excluding $\set{0}$
if $p<0$), i.e. $E_p(f)(t\theta) = t^p f(\theta)$ for $t>0$ and
$\theta \in S^{n-1}$. We denote by $E_p^\wedge(f)$ the Fourier
transform of $E_p(f)$ as a distribution. Note that $E_p^\wedge(f)$
need not necessarily be a continuous function on $\Real^n \setminus
\set{0}$, nor even a measure on $\Real^n$. In order to ensure that
$E_p^\wedge(f)$ is a continuous function, we need to add some
smoothness assumptions on $f$ (\cite{Koldobsky-Book}). We remark
that for an infinitely smooth function $f \in C(S^{n-1})$,
$E_p^\wedge(f)$ is infinitely smooth on $\Real^n \setminus \set{0}$
for any $p\in(-n,0)$. Whenever $E_p^\wedge(f)$ is continuous on
$\Real^n \setminus \set{0}$, it is uniquely determined by its value
on $S^{n-1}$ (by homogeneity), so we identify (abusing notation)
between $E_p^\wedge(f)$ and its restriction to $S^{n-1}$.

Clearly $E_{-k}(\rho_K^k) = \norm{\cdot}_K^{-k}$ for a star-body
$K$. Given a full-rank linear transformation $T$ in $\Real^n$, we
denote $T(E_p(f)) = E_p(f) \circ T^{-1}$, so $T(E_{-k}(\rho_K^k) =
E_{-k}(\rho_{T(K)}^k)$ for a star-body $K$. Again, we identify (by
homogeneity) between $T(E_p(f))$ and its restriction on $S^{n-1}$.

It is easy to check (e.g.
\cite{EMilman-Generalized-Intersection-Bodies}) that for any
infinitely smooth $K \in D_{-k}^n$, we have $E_{-k}^\wedge(\rho_K^k)
\geq 0$ (and clearly $\rho_K^k \geq 0$). In fact, this immediately
follows from the fact that this is true for $D_n \in D_{-k}^n$, the
linearity of the Fourier transform, and its behavior under full-rank
linear transformations. With Theorem \ref{thm:G&Z} and
characterization (2) of Theorem \ref{thm:I-char} in mind, asking
whether $\BP_k^n = \I_k^n$ is equivalent to asking whether the only
infinitely smooth functions $f \in C(S^{n-1})$ such that $f \geq 0$
and $E_{-k}^\wedge(f) \geq 0$, are the ones such that $f = \rho_K^k$
for some $K \in D_{-k}^n$. In other words, whether every such $f$
can be approximated (in the maximum norm in $C(S^{n-1})$, which is
clearly the same for $f$ and for $f^{1/k}$) by functions of the form
$\sum_{i=1}^m T_i(E_{-k}(1))$, where $T_i$ are full-rank linear
transformations. The following is thus an immediate consequence of
Theorem \ref{thm:main1}:

\begin{cor}
Let $n \geq 4$ and $2 \leq k \leq n-2$. Then there exists a
``non-trivial" infinitely smooth function of revolution $f \in
C(S^{n-1})$ such that $f \geq 0$ and $E_{-k}^\wedge(f) \geq 0$. By
``non-trivial", we mean that $f$ cannot be approximated in the
maximum norm on $C(S^{n-1})$ by functions of the form $\sum_{i=1}^m
T_i(E_{-k}(1))$, where $\set{T_i}$ are full-rank linear
transformations in $\Real^n$.
\end{cor}

\medskip

To conclude, we comment that although the original definitions of
$\BP_k^n$ and $\I_k^n$ make sense only for integer values of $k$
(between $1$ and $n-1$), some of the alternative characterizations
of these classes stated in this section make sense for arbitrary
real-valued $k$, for $0<k<n$. In particular, characterizations (2)
and (3) of Theorem \ref{thm:I-char} for the class $\I_k^n$ and
Theorem \ref{thm:G&Z} for the class $\BP_k^n$ may be taken as
definitions for these classes of star-bodies in this extended range
of $k$. It then makes sense to ask whether Theorem \ref{thm:main1}
also holds for any non-integer $1<k<n-1$. Although we do not proceed
in this direction, the answer should be positive, since our
construction of the function $\tilde{g}$ in Proposition
\ref{prop:main} is purely analytic, and everything still works for
arbitrary real-valued $k$, for $1<k<n$.


\bibliographystyle{amsalpha}
\bibliography{../../ConvexBib}

\end{document}